\documentclass[12pt]{article}
\usepackage{amsthm,amsfonts,amssymb,amscd}

\begin{document}

\begin{center}
\Large \bf On the multiplicity of solutions \\
of a system of algebraic equations
\end{center}
\vspace{1cm}

\centerline{A.V.Pukhlikov
}

\parshape=1
3cm 10cm \noindent {\small \quad \quad \quad
\quad\quad\quad\quad\quad\quad\quad {\bf }\newline We obtain upper
bounds for the multiplicity of an isolated solution of a system of
equations $f_1=\dots = f_M =0$ in $M$ variables, where the set of
polynomials $(f_1,\dots, f_M)$ is a tuple of general position in a
subvariety of a given codimension which does not exceed $M$, in
the space of tuples of polynomials. It is proved that for
$M\to\infty$ that multiplicity grows not faster than
$\sqrt{M}\exp[\omega\sqrt{M}]$, where $\omega>0$ is a certain
constant.

Bibliography: 3 titles.} \vspace{1cm}

{\bf Introduction}\vspace{0.5cm}

In the present paper, the following problem is considered. Let
$$
\left\{\begin{array}{c} f_1(z_1,\dots,z_M)=0\\ \dots \\
f_M(z_1,\dots,z_M)=0
\end{array}\right.
$$
be a system of polynomial equations of degree $d\geq2$, which has
the origin $o=(0,\dots, 0)\in{\mathbb C}^M$ as an isolated
solution. For a given $m\geq 1$ one needs to estimate the
codimension of the set of such tuples $(f_1,\dots,f_M)$, that
$$
\mathop{\rm dim} {\cal O}_{o,{\mathbb C}^M}/(f_1,\dots,f_M)\geq m,
$$
in the space of all tuples of polynomials of degree $d$ with no
free term. Informally speaking, how many independent conditions on
the coefficients of the polynomials $f_1,\dots,f_M$ are imposed if
it is required that the multiplicity of the given solution is no
smaller than $m$? Problems of that type emerge in the theory of
birational rigidity (see [1, Proposition 3.3]). As another
application, we point out the problem of description of possible
singularities of the variety of lines on a generic Fano variety
$V\subset {\mathbb P}^N$ in a given family. However, this problem
is interesting by itself, too. The problem described above can be
formulated in another way: for a given codimension $a\geq 1$ to
estimate the maximal possible multiplicity for a {\it generic}
tuple of equations $(f_1,\dots,f_M)\in B$ in a {\it given}
subvariety $B$ of codimension $a$. Thus we are looking for the
maximum over all subvarieties $B$, in each of which a tuple of
general position is taken. Of course, this problem makes sense
only provided that the set of such tuples $(f_1,\dots,f_M)$, that
the system $f_1=\dots=f_M=0$ has a set of solutions of positive
dimension, containing the point $o$, is of codimension not less
than $a+1$. This is true, if $a\leq M$.\vspace{0.1cm}

In [1,\S3] a simple example (the idea of which is actively used in
this paper) was constructed, which shows that for $M\gg 0$ the
maximal multiplicity of an isolated solution in codimension $a=M$
grows not slower than
$$
2^{\sqrt{M}}.
$$
In the present paper for this value we obtain the upper bound
$$
\sqrt{M}e^{\omega\sqrt{M}},
$$
where $\omega>0$ is a certain concrete real number. To do this, we
generalize the problem above for systems of $i\leq M$ polynomial
equations, which makes it possible to construct an inductive
procedure of estimating the maximal intersection multiplicity for
a given codimension of the set of equations.\vspace{0.1cm}

Let us explain the main difficulty in solving the problem above.
Let $Y_i$, $i=1,\dots,M$, be the algebraic cycle of the
scheme-theoretic intersection
$$
(\{f_1=0\}\circ\dots\circ\{f_i=0\})
$$
in a neighborhood of the point $o$. This is an effective cycle of
codimension $i$. Set $m_i=\mathop{\rm mult}\nolimits_{o} Y_i$. It
seems natural to consider the whole sequence of multiplicities
$(m_1,\dots,m_M)$, estimating the codimension of the space of
polynomials $f_{i+1}$ in terms of the jump of the multiplicity
from $m_i$ to $m_{i+1}$ (this very approach was realized in
[1,\S3]). However, in our problem this approach does not
work.\vspace{0.1cm}

Let $\widetilde{{\mathbb C}^M}\to {\mathbb C}^M$ be the blow up of
the point $o$, $E\cong {\mathbb P}^{M-1}$ the exceptional divisor,
$\widetilde{Y}_i$ the strict transform of $Y_i$,
$(\widetilde{Y}_i\circ E)=\sum c_jR_j$ the algebraic projectivized
tangent cone. According to the intersection theory [2], the
multiplicity of the scheme-theoretic intersection of the cycle
$Y_i$ and the divisor $D_{i+1}=\{f_{i+1}=0\}$ at the point $o$ is
given by the formula
$$
m_{i+1}=m_i\mathop{\rm
mult}\nolimits_{o}D_{i+1}+\sum_{R_{jk}}d_{jk}\left(\mathop{\rm
mult}\nolimits_{R_{jk}}\widetilde{Y}\right)\left(\mathop{\rm
mult}\nolimits_{R_{jk}}\widetilde{D_{i+1}}\right),
$$
where the sum is taken over some finite set of irreducible
subvarieties of codimension $(i+1)$, including infinitely near
ones, $R_{jk}$ covers $R_j$ with the multiplicity $d_{jk}$. Taking
into account that $M\gg 0$, for $i$ close to $M$ the structure of
the singularity of the cycle $Y_i$ at the point $o$ can be
arbitrary, that is, it can not be explicitly described. Therefore,
it is impossible to estimate, how many independent conditions on
the polynomial $f_{i+1}$ for $f_1,\dots,f_i$ fixed are imposed by
the bounds for the multiplicities $\mathop{\rm
mult}\nolimits_{R_{jk}}\widetilde{D_{i+1}}$. The only and obvious
conclusion, which can be derived from the formula for $m_{i+1}$,
given above, is that the condition $m_{i+1}\geq c$ for a fixed
cycle $Y_i$ defines a closed subset in the space of polynomials
$f_{i+1}$, which is a union of a finite number of linear
subspaces. Indeed, the condition $\mathop{\rm
mult}\nolimits_{R_{jk}}\widetilde{D_{i+1}}\geq\gamma$ is a linear
one.

By what was said above, in order to get an effective bound for the
maximal intersection multiplicity in codimension $a\geq 1$ one
needs a different approach, which is developed in the present
paper. The main idea is to estimate the maximal multiplicity for
$i$ polynomials via the maximal multiplicity for $(i-1)$
polynomials with an appropriate correction of the codimension. The
estimates, obtained by means of this inductive method, seem to be
close to the optimal ones.\vspace{0.1cm}

The paper is organized in the following way. In \S1 we develop an
inductive procedure of estimating the multiplicity. Using it, in
\S2 we derive an absolute estimate of the intersection
multiplicity and, as a corollary, the main asymptotic result of
this paper. In \S3, following [3], we briefly remind the method of
estimating the codimension of the set of tuples $(f_1,\dots,f_i)$,
defining sets of an ``incorrect'' codimen\-sion $\leq
i-1$.\vspace{0.1cm}

To conclude, we note that the problem, considered in this paper,
can be set up and solved by the same method for an arbitrary very
ample class $H$ on an algebraic variety $V$ at a point $o\in
V$.\vspace{1cm}

{\bf \S1. The inductive method of estimating the
multiplicity}\vspace{0.5cm}

In this section we develop an inductive procedure of estimating
the maximal multiplicity in a given codimension. In the beginning
of the section we consider equations of arbitrary degree $d\geq
2$, later we restrict ourselves by quadratic polynomi\-als
$(d=2)$. For a codimension, not exceeding $M$, this does not
change the result (see Remark 1.4).\vspace{0.3cm}

{\bf 1.1. Set up of the problem.} Fix the complex coordinate space
${\mathbb C}^M_{(z_1,\dots,z_M)}$, $M\geq 1$. By the symbol ${\cal
P}_{d,M}$ we denote the space of {\it homogeneous} polynomials of
degree $d\geq1$ in the variables $z_*$, by the symbol ${\cal
P}_{\leq d,M}$ we denote the space of polynomials of degree $\leq
d$ {\it with no free term} in the variables $z_*$. On each of
these spaces there is a natural action of the matrix group
$GL_M({\mathbb C})$. Set
$$
{\cal P}^{i}_{\leq d,M}=\underbrace{{\cal P}_{\leq d,M}\times\dots
\times {\cal P}_{\leq d,M}}_i
$$
to be the space of tuples $(f_1,\dots,f_i)$. By the symbol
$$
Z(f_1,\dots,f_i)
$$
we denote the subscheme $\{f_1=\dots=f_i=0\}$, which we will study
in a neighborhood of the point $o=(0,\dots,0)$, that is, in fact,
the subject of our study is the local ring
$$
{\cal O}_{o,Z(f_*)}={\cal O}_{o,{\mathbb C}^M}/(f_1,\dots,f_i).
$$
Denote the map
$$
\mu\colon {\cal P}^{i}_{\leq d,M}\to {\mathbb Z}_+\cup\{\infty\},
$$
setting $\mu(f_1,\dots,f_i)=\infty$, if $\mathop{\rm
codim}\nolimits_o Z(f_1,\dots,f_i)\leq i-1$ (the symbol
$\mathop{\rm codim}\nolimits_o$ stands for the codimension in a
neighborhood of the point $o$), and
$$
\mu(f_1,\dots,f_i)=\mathop{\rm mult}\nolimits_{o}
Z(f_1,\dots,f_i),
$$
if $\mathop{\rm codim}\nolimits_o Z(f_1,\dots,f_i)=i$. For an
arbitrary irreducible subvariety $B\subset {\cal P}^{i}_{\leq
d,M}$ set
$$
\mu(B)=\mathop{\rm min}\limits_{(f_1,\dots,f_i)\in B}
\{\mu(f_1,\dots,f_i)\}\in {\mathbb Z}_+\cup\{\infty\}.
$$
Therefore, $\mu(B)=\infty$ if and only if for every tuple of
polynomials $(f_1,\dots,f_i)\in B$ the complete intersection
$Z(f_1,\dots,f_i)$ has in a neighborhood of the point $o$ an
``incorrect'' codimension $\leq i-1$. The equality
$\mu(B)=m\in{\mathbb Z}_+$ means that for a generic tuple of
polynomials $(f_1,\dots,f_i)\in B$ the complete intersection
$Z(f_1,\dots,f_i)$ has in a neighborhood of the point $o$ the
correct codimension $i$ and its multiplicity at the point $o$ is
$m\geq 1$.\vspace{0.1cm}

{\bf Definition 1.1.} The maximal intersection multiplicity of a
generic tuple of polynomials at the point $o$ {\it in the
codimension} $a\in{\mathbb Z}_+$ is
$$
\mu_i(a)=\mathop{\rm max}\limits_{B\subset{\cal P}^{i}_{\leq d,M}}
\{\mu(B)\}\in {\mathbb Z}_+\cup\{\infty\},
$$
where the maximum is taken over all irreducible subvarieties
$B\subset{\cal P}^{i}_{\leq d,M}$ of codimension
$a$.\vspace{0.1cm}

Definition 1.1 can be re-formulated as follows. The multiplicity
$\mu_i(a)$ is $\infty$, if and only if the codimension of the
closed algebraic set
$$
\{(f_1,\dots,f_i)\in{\cal P}^{i}_{\leq
d,M}\,|\,\mu(f_1,\dots,f_i)=\infty\}
$$
does not exceed $a$ (and in that case for $B$ we can take any
irreducible subvariety of codimension $a$, contained in that set).
Otherwise, the multiplicity $\mu_i(a)$ is the minimal positive
integer $m\geq 1$, satisfying the condition: the codimension of
the closed algebraic set
$$
\{(f_1,\dots,f_i)\in{\cal P}^{i}_{\leq
d,M}\,|\,\mu(f_1,\dots,f_i)\geq m+1\}
$$
is not less than $a+1$. In other words, for any irreducible
subvariety $B\subset {\cal P}^{i}_{\leq d,M}$ of codimension $a$
and a generic tuple $(f_1,\dots,f_i)\in B$ we get
$$
\mu(f_1,\dots,f_i)\leq\mu_i(a)
$$
and for a certain subvariety $B$ this inequality turns into the
equality.\vspace{0.1cm}

{\bf Remark 1.1.} Apart from the matrix group $GL_M({\mathbb C})$,
which acts naturally on the space ${\cal P}^{i}_{\leq d,M}$ by
linear changes of coordinates, on that space naturally acts the
matrix group $GL_i({\mathbb C})$: with a non-degenerate $(i\times
i)$ matrix $A$ we associate the transformation of the tuple of
polynomials
$$
(f_1,\dots,f_i)\mapsto(f_1,\dots,f_i)A.
$$
The multiplicity $\mu(f_1,\dots,f_i)$ is invariant with respect to
the action of these two groups. Respectively, the algebraic sets
$$
X_{i,M}(m)=\{(f_1,\dots,f_i)\,|\, \mu(f_1,\dots,f_i)\geq
m\}\subset {\cal P}^{i}_{\leq d,M}
$$
and their irreducible components are $GL_M({\mathbb C})$- and
$GL_i({\mathbb C})$-invariant. For this reason, the definition of
the number $\mu_i(a)$ can be modified in the following way: for
any $GL_M({\mathbb C})$- and $GL_i({\mathbb C})$-invariant
subvariety $B\subset {\cal P}^{i}_{\leq d,M}$ of codimension $\leq
a$ we have $\mu(B)\leq\mu_i(a)$, and moreover, for a certain
(invariant) $B$ this is an equality. The equivalence of the two
definitions of the number $\mu_i(a)$ is obvious: let us consider
the closed set $X_{i,M}(\mu_i(a))$. Its codimension in the space
${\cal P}^{i}_{\leq d,M}$ does not exceed $a$ and each of its
components is invariant, and moreover, for some component $B$ of
codimension $\leq a$ we have $\mu(B)=\mu_i(a)$, which is what we
need.

Now let us consider the problem, for which values $a\in{\mathbb
Z}_+$ the numbers $\mu_i(a)$ are certainly finite.\vspace{0.1cm}

{\bf Proposition 1.1.} {\it The codimension of the closed set
$X_{i,M}(\infty)$ for $i\leq M-1$ is not less than $dM$, and for
$i=M$ not less than $(d-1)M+1$}.\vspace{0.1cm}

{\bf Proof} is given in \S3.\vspace{0.1cm}

{\bf Corollary 1.1.} {\it For $a\leq M$ we have
$\mu_i(a)<\infty$}.\vspace{0.1cm}

The problem of estimating the numbers $\mu_i(a)$ from above is
considered in this paper for those values of $a$
only.\vspace{0.3cm}


{\bf 1.2. The invariant $\varepsilon$ and reduction to the
standard form.} For an irreducible subvariety $B\subset {\cal
P}^{i}_{\leq 2,M}$ we define the number
$$
\varepsilon(B)=i -\mathop{\rm
rk}(df_1(o),\dots,df_i(o))\in\{0,1,\dots,i\},
$$
where $(f_1,\dots,f_i)\in B$ is a tuple of general position. If
the subvariety $B$ is $GL_i({\mathbb C})$-invariant, then the
equality $\varepsilon(B)=b$ means that in a generic tuple
$(f_1,\dots,f_i)\in B$ the first $(i-b)$ linear forms
$$
df_1(o),\dots,df_{i-b}(o)
$$
are linearly independent, whereas the forms $df_{i-b+j}(o)$ for
$j\in\{1,\dots,b\}$ are their linear combinations. For any
irreducible subvariety $B$, satisfying the latter condition, there
exists a non-empty Zariski open subset $B^{o}\subset B$, on which
the map of {\it reducing to the standard form} is well defined:
$$
\rho\colon B^{o}\to {\cal P}^{i-b}_{\leq 2,M}\times {\cal
P}^{b}_{2,M},
$$
$$
\rho\colon (f_1,\dots,f_i)\mapsto
(f_1,\dots,f_{i-b},f^+_{i-b+1},\dots,f^+_i),
$$
where
$$
f^+_{i-b+j}=f_{i-b+j}-\sum^{i-b}_{\alpha=1}\lambda_{j\alpha}f_{\alpha},
$$
the coefficients $\lambda_{j\alpha}$ are defined by the equalities
$$
df_{i-b+j}(o)=\sum^{i-b}_{\alpha=1}\lambda_{j\alpha}df_{\alpha}(o).
$$
Therefore, $df^+_{i-b+j}(o)=0$ and $f^+_{i-b+j}$ is a homogeneous
polynomial of degree 2. The closure of the image $\rho(B^o)$ we
denote by the symbol $\bar B$. If the subvariety $B$ is invariant
with respect to the action of $GL_i({\mathbb C})$, then the
coefficients $\lambda_{j\alpha}$ take arbitrary values and for
that reason
$$
\mathop{\rm dim}B=\mathop{\rm dim}\bar{B}+b(i-b).
$$
Obviously, every fibre of the map $\rho\colon B^o\mapsto\rho(B^o)$
is ${\mathbb C}^{b(i-b)}$.

The main technical tool for estimating the numbers $\mu_i(a)$ is
given by the more sensitive numbers
$$
\mu_{i,M}(a,b)=\mathop{\rm max}\left\{ m\geq 1\, \left| \,
\begin{array}{l}
\mbox{the set}\,\, X_{i,M}\,\, \mbox{has an irreducible component}
\\
B \,\,\mbox{of codimension}\,\,\leq a
\,\,\mbox{with}\,\,\varepsilon(B)=b
\end{array}\right.\right\}.
$$
Obviously, $\mu_i(a)=\mathop{\rm max}\limits_b\{\mu_{i,M}(a,b)\}$,
where the maximum is taken over all possible values of the number
$\varepsilon(B)$ for irreducible subvarieties $B$ of codimension
$\leq a$. It is easy to see that the codimension of the subset
$$
\{(f_1,\dots,f_i)\,|\,\mathop{\rm rk}(df_1(o),\dots,df_i(o))\leq
i-b\}
$$
is $b(M+b-i)$, so that the equality $\varepsilon(B)=b$ is only
possible if $a\geq b(M+b-i)$. In the sequel, when the notation
$\mu_{i,M}(a,b)$ is used, it means automatically that the latter
inequality holds. The following obvious fact is
true.\vspace{0.1cm}

{\bf Proposition 1.2.} {\it The equality
$$
\mu_{i,M}(a,0)=1
$$
holds.}

{\bf Proof.} If $\varepsilon(B)=0$, then for a generic tuple
$(f_1,\dots,f_i)$ the differentials $df_1(o),\dots,df_i(o)$ are
linearly independent, that is, the set $\{f_1=\dots=f_i=0\}$ is a
smooth subvariety of codimension $i$ in a neighborhood of the
point $o$, which is what we need. Q.E.D.

Let us find an upper bound for the numbers $\mu_{i,M}(a,b)$ for
$b\geq 1$.\vspace{0.3cm}


{\bf 1.3. Splitting off the last factor.} Let
$$
\pi_i\colon {\cal P}^{i-b}_{\leq 2,M}\times {\cal P}^{b}_{2,M} \to
{\cal P}^{i-b}_{\leq 2,M}\times {\cal P}^{b-1}_{2,M}
$$
be the projection along the last direct factor ${\cal P}_{2,M}$.
For the closed set $\bar{B}\subset {\cal P}^{i-b}_{\leq 2,M}\times
{\cal P}^{b}_{2,M}$, constructed above, denote by the symbol
$[\bar B]_{i-1}$ the closure of the set $\pi_i(\bar B)$. It is
easy to see that the following relation holds:
$$
\mathop{\rm codim}\left(\bar{B}\subset {\cal P}^{i-b}_{\leq
2,M}\times {\cal P}^{b}_{2,M}\right)= \mathop{\rm codim}\left(
B\subset {\cal P}^{i}_{\leq 2,M}\right)-(M+b-i)b.
$$
Starting from this moment, unless otherwise specified, the
codimension is always meant with respect to the natural ambient
space; for instance, the last equality writes simply as
$$
\mathop{\rm codim}\bar{B}=\mathop{\rm codim} B-(M+b-i)b.
$$
Sometimes for the convenience of the reader we remind, the
codimension with respect to which space is meant.

For a tuple of general position $(f_1,\dots,f_{i-1})\in[\bar
B]_{i-1}$ denote by the symbol
$$
[\bar{B}]^i=[\bar{B}]^i(f_1,\dots,f_{i-1})\subset {\cal P}_{2,M}
$$
the fibre of the projection $\pi_i|_{\bar B}\colon \bar
B\mapsto[\bar B]_{i-1}$. Obviously,
$$
\mathop{\rm codim}\bar{B}=\mathop{\rm codim}[\bar{B}]_{i-1}+
\mathop{\rm codim}[\bar{B}]^i
$$
(recall: the codimension is meant with respect to the natural
ambient space, for $\bar B$ it is ${\cal P}^{i-b}_{\leq 2,M}\times
{\cal P}^{b}_{2,M}$, for $[\bar B]_{i-1}$ it is the space ${\cal
P}^{i-b}_{\leq 2,M}\times {\cal P}^{b-1}_{2,M}$, for $[\bar B]^i$
it is ${\cal P}_{2,M})$. Set
$$
\gamma_i=\gamma_i(B)=\mathop{\rm codim}[\bar{B}]^i.
$$
Since $\mathop{\rm codim} B\leq a$, we obtain the estimate
$$
\mathop{\rm codim}[\bar{B}]_{i-1}=\mathop{\rm codim}
B-(M+b-i)b-\gamma_i\leq
$$
$$
\leq a-(M+b-i)b-\gamma_i.
$$
This, in particular, implies that
$$
0\leq \gamma_i\leq a-(M+b-i)b.
$$
\vspace{0.3cm}

\newpage

{\bf 1.4. The main inductive estimate.} The following fact is
true.\vspace{0.1cm}

{\bf Theorem 1.} {\it For any $i,M,a,b$ there exist integers
$\alpha\in\{0,1\}$ and $\gamma\in\{0,\dots,a-(M+b-i)b\}$ such that
the following inequality holds:}
\begin{equation}\label{june2011.1}
\begin{array}{ccl}
\mu_{i,M}(a,b) & \leq & \mu_{i-1,M}(a-(M+b-i)-\gamma,b-1)+\\ & &
\\ & &
+\mu_{i-1,M-1}(a-(M+b-i)-\alpha(b-1),b-\alpha).\end{array}
\end{equation}

{\bf Remark 1.2.} As we will see from the proof of the theorem,
the numbers $\alpha$ and $\gamma$ are determined by the subvariety
$B$, which realizes the multiplicity $\mu_{i,M}(a,b)$. There can
be more than one such subvariety; respectively, several
inequalities (\ref{june2011.1}) can be satisfied for the number
$\mu_{i,M}(a,b)$, with different values of $\alpha$ and $\gamma$.
Furthermore, the inequalities
$$
\mu_{i,M}(a_1,b)\leq
\mu_{i,M}(a_2,b)\quad\mbox{and}\quad\mu_i(a_1)\leq\mu_i(a_2)
$$
hold for $a_1\leq a_2$, which implies that in (\ref{june2011.1})
one can set $\gamma=0$ and the estimate still holds (possibly
becomes weaker).\vspace{0.1cm}

{\bf Proof of Theorem 1.} Let us fix a $GL_i({\mathbb
C})$-invariant irreducible subvariety $B$, realizing the value
$\mu_{i,M}(a,b)$, $\varepsilon(B)=b$. We may assume that $B$ is an
irreducible component of the closed set $X_{i,M}(m)$, where
$m=\mu_{i,M}(a,b)$. To simplify the formulas, we assume that
$\mathop{\rm codim}B=a$ (if $\mathop{\rm codim}B<a$, then the
estimates below can only become stronger). Fix a linear form
$L(z_1,\dots,z_M)$ of general position. In particular, if
$(f_1,\dots,f_i)\in B$ is a generic tuple, so that the set
$\{f_1=\dots=f_{i-1}=0\}$ is of codimension $(i-1)$ in a
neighborhood of the point $o$, the multiplicity of the effective
cycle
$$
(\{f_1=0\}\circ\dots\circ\{f_{i-1}=0\})
$$
at the point $o$ is equal to the multiplicity of the intersection
of that cycle with the hyperplane $\{L=0\}$ at the point $o$. Let
$$
\Pi_L=\{L(z_*)L_1(z_*)\,|\, L_1\in {\cal P}_{1,M}\}\subset {\cal
P}_{2,M}
$$
be the linear space of reducible homogeneous quadratic
polynomials, divisible by $L$. Set
$$
{\cal P}_L={\cal P}^{i-b}_{\leq 2,M}\times {\cal
P}^{b-1}_{2,M}\times \Pi_L\subset {\cal P}^{i-b}_{\leq 2,M}\times
{\cal P}^{b}_{2,M}.
$$
This is a closed subset. The intersection $\bar{B}\cap {\cal P}_L$
is non-empty and of codimension not higher than $\mathop{\rm
codim}\bar{B}$ in ${\cal P}_L$. By the symbol $[\bar{B}\cap {\cal
P}_L]_{i-1}$ we denote the closure of the set $\pi_i(\bar{B}\cap
{\cal P}_L)$. As we consider only codimensions $a\leq M$, the
equality
$$
[\bar{B}\cap {\cal P}_L]_{i-1}=[\bar{B}]_{i-1}
$$
holds, since for a generic tuple $(f_1,\dots,f_i)\in B$ the
intersection of the space $\Pi_L$ with the fibre $[\bar
B]^i(f_1,\dots,f_{i-1})$ has a positive dimension. More precisely,
the codimension of that intersection in $\Pi_L\cong {\cal
P}_{1,M}$ does not exceed $\gamma_i$.\vspace{0.1cm}

{\bf Remark 1.3.} Since we assume that $B$ is an irreducible
component of the closed set $X_{i,M}(m)$, the fibre
$$
[\bar{B}]^i=\{f_i\in {\cal P}_{2,M}\, |\, \mathop{\rm
mult}\nolimits_{o}\{f_1=\dots =f_i=0\}\geq m\},
$$
$m=\mu_{i,M}(a,b)$, for $f_1,\dots,f_{i-1}$ fixed, is a union of a
finite number of linear subspaces of codimension $\gamma_i$.
Therefore, the closed set
$$
\Pi(f_1,\dots,f_{i-1})=\{L_1\in {\cal P}_{1,M}\, |\, \mathop{\rm
mult}\nolimits_{o}\{f_1=\dots =f_{i-1}=LL_1=0\}\geq m\}
$$
for a generic tuple $(f_1,\dots,f_{i-1})\in[\bar B]_{i-1}$ is a
union of a finite number of linear subspaces in ${\cal P}_{1,M}$,
the codimension of each of which in ${\cal P}_{1,M}$ does not
exceed $\gamma_i$.\vspace{0.1cm}

By what was said, the inequality
$$
\begin{array}{rcl}
m=\mu_{i,M}(a,b) & \leq & \mathop{\rm
mult}\nolimits_{o}\{f_1=\dots =f_{i-1}=L=0\}+ \\ & & \\
   &  +  & \mathop{\rm
mult}\nolimits_{o}\{f_1=\dots =f_{i-1}=L_1=0\}
\end{array}
$$
holds. Since $L$ is a form of general position, the first summand
in the right hand side is
$$
\mathop{\rm mult}\nolimits_o\{f_1=\dots=f_{i-1}=0\}.
$$
Here $(f_1,\dots,f_{i-1})\in[\bar B]_{i-1}$ is a tuple of general
position. Now the set $[\bar{B}]_{i-1}\subset {\cal P}^{i-b}_{\leq
2,M}\times {\cal P}^{b-1}_{2,M}$ can be represented as a result of
reducing to the standard form of the closed subset $C\subset {\cal
P}^{i-1}_{\leq 2,M}$, that is,
$$
[\bar B]_{i-1}=\bar C,
$$
where $C$ is constructed by the procedure, which is inverse to the
procedure of reducing to the standard form: $C$ is the closure of
the set of $(i-1)$-tuples
$$
\{(g_1,\dots,g_{i-b},g^+_{i-b+1},\dots,g^+_{i-1})\},
$$
where
$$
g^+_{i-b+j}=g_{i-b+j}+\sum^{i-b}_{\alpha=1}\lambda_{j\alpha}g_{\alpha},
$$
for all $(g_1,\dots,g_{i-1})\in[\bar B]_{i-1}$ and
$\lambda_{j\alpha}\in {\mathbb C}$. From this, it follows that
$$
\mathop{\rm dim} C=\mathop{\rm dim} [\bar{B}]_{i-1}+(b-1)(i-b),
$$
so that
$$
\mathop{\rm codim} C = \mathop{\rm codim} [\bar{B}]_{i-1}
+(b-1)(M+b-i)=
$$
$$
=\mathop{\rm codim} B - (M+b-i)-\gamma_i.
$$
(Recall, that each of the three codimensions is taken with respect
of the {\it corresponding} ambient space; for instance, for $C$ it
is ${\cal P}^{i-1}_{\leq 2,M}$). Since, obviously,
$\varepsilon(C)=b-1$, we obtain that
$$
\mathop{\rm mult}\nolimits_{o}\{f_1=\dots =f_{i-1}=0\}\leq
\mu_{i-1,M}(a-(M+b-i)-\gamma_i, b-1).
$$
This gives us the first half of the right hand side of the
inequality of Theorem 1.\vspace{0.3cm}

{\bf 1.5. The multiplicity of intersection with the hyperplane
$\{L_1=0\}$.} It remains to estimate the multiplicity $\mathop{\rm
mult}\nolimits_{o}\{f_1,\dots =f_{i-1}=L_1=0\}$. This is somewhat
harder, since the form $L_1$ depends on the tuple
$(f_1,\dots,f_{i-1})$ and for this reason is not a form of general
position with respect to that tuple. Note that for a generic tuple
$(f_1,\dots,f_{i-1})$ the set $\Pi(f_1,\dots,f_{i-1})$ does not
depend on the choice of the form $L$. Therefore, the set
$$
\Pi\subset {\cal P}^{i-b}_{\leq 2,M}\times {\cal
P}^{b-1}_{2,M}\times {\cal P}_{1,M},
$$
defined as the closure of the set of tuples
$$
(f_1,\dots,f_{i-1},L_1\in\Pi(f_1,\dots,f_{i-1}))
$$
for generic tuples $(f_1,\dots,f_{i-1})\in[\bar{B}]_{i-1}$, does
not depend on the choice of the form $L$, either. Since that form
of general position $L$ does not take part in the subsequent
constructions, to simplify the notations we write $L$ instead of
$L_1$, if it does not generate a confusion.

Obviously, the set $\Pi$ is invariant with respect to the action
of the group $GL_M({\mathbb C})$, therefore the projection
$$
\pi\colon \Pi\to {\cal P}_{1,M},
$$
$$
\pi\colon(f_1,\dots,f_{i-1},L)\mapsto L,
$$
is surjective and all its fibres are of the same dimension. Since
the codimension of the closed set $\Pi$ (with respect to the
ambient space ${\cal P}^{i-b}_{\leq 2,M}\times {\cal
P}^{b-1}_{2,M}\times {\cal P}_{1,M}$) does not exceed the number
$$
\mathop{\rm codim}[\bar{B}]_{i-1}+\gamma_i=a-(M+b-i)b,
$$
for a generic linear form $L\in {\cal P}_{1,M}$ the codimension of
the fibre $\pi^{-1}(L)\subset {\cal P}^{i-b}_{\leq 2,M}\times
{\cal P}^{b-1}_{2,M}$ is bounded from above by the same number
$a-(M+b-i)b$.

Now for a generic tuple $(f_1,\dots,f_{i-1})\in\pi^{-1}(L)$ there
are two options:\vspace{0.1cm}

1) either the differentials
$(df_1|_{\{L=0\}}(o),\dots,df_{i-b}|_{\{L=0\}}(o))$ remain
linearly independent (an equivalent formulation: the subspace
\begin{equation}\label{june2011.2}
\{df_1(o)=\dots=df_{i-b}(o)=0\}
\end{equation}
is not contained in the hyperplane $\{L=0\}$),\vspace{0.1cm}

2) or the rank of the set of linear forms
$$
df_1(o)|_{\{L=0\}},\dots,df_{i-b}(o)|_{\{L=0\}}
$$
drops by one (an equivalent formulation: the subspace
(\ref{june2011.2}) is contained in the hyperplane
$\{L=0\}$).\vspace{0.1cm}

In the case 1) set $\alpha=\alpha(B)=1$, in the case 2) set
$\alpha=\alpha(B)=0$. Furthermore, let
$$
\bar{B}_L\subset {\cal P}^{i-b}_{\leq 2,M-1}\times {\cal
P}^{b-1}_{2,M-1}
$$
be the closure of the set
$$
\{(f_1|_{\{L=0\}},\dots,f_{i-1}|_{\{L=0\}})\,|\,\,
(f_1,\dots,f_{i-1})\in\pi^{-1}(L)\}.
$$

Let us consider first the case 1). Here for a generic tuple
$(g_1,\dots,g_{i-1})\in{\bar B}_L$ the differentials
$$
dg_1(o),\dots,dg_{i-b}(o)
$$
are linearly independent, and for $j\geq i-b+1$ we have
$dg_j(o)=0$. Now we argue as in Sec. 1.4: the set ${\bar B}_L$ is
the result of reducing to the standard form of a certain closed
set $C\subset {\cal P}^{i-1}_{\leq 2,M-1}$. The set $C$ is
obtained from ${\bar B}_L$ by the procedure, which is converse to
the procedure of reducing to the standard form. Obviously,
$\varepsilon(C)=b-1$ and
$$
\mathop{\rm codim} C = \mathop{\rm codim}\bar{B}_L- (i-b)(b-1)+
(M-1)(b-1)=
$$
$$
= \mathop{\rm codim} \bar{B}_L + (M+b-i-1)(b-1),
$$
so that taking into account the estimate
$$
\mathop{\rm codim} \bar{B}_L\leq \mathop{\rm codim}
\pi^{-1}(L)\leq a-(M+b-i)b
$$
we obtain the inequality
$$
\mathop{\rm codim} C\leq a-(M+b-i)-\alpha (b-1).
$$
Since
$$
\mathop{\rm mult}\nolimits_{o}\{f_1=\dots=f_{i-1}=L=0\}=
\mathop{\rm
mult}\nolimits_{o}\{f_1|_{\{L=0\}}=\dots=f_{i-1}|_{\{L=0\}}=0\},
$$
we obtain the final upper estimate for that multiplicity: it can
not exceed the number
$$
\mu_{i-1,M-1}(a-(M+b-i)-\alpha(b-1),e-\alpha).
$$
(Recall that in the case under consideration $\alpha=1$, and in
the inequalities above the codimension is taken with respect to
the natural ambient spaces, each of the sets ${\bar B}_L$, $C$,
$\pi^{-1}(L)$ has its own ambient space.)\vspace{0.1cm}

Now let us consider the case 2). Here for a generic tuple
$(g_1,\dots,g_{i-1})\in{\bar B}_L$ the rank of the system of
linear functions $dg_1(o),\dots,dg_{i-b}(o)$ is equal to $i-b-1$.
We may assume that the first $i-b-1$ of them are linearly
independent, and $dg_{i-b}(o)$ is their linear combination. For
$j\geq i-b+1$ we get, as above, that $dg_j(o)=0$. In the case 2)
the set ${\bar B}_L$ is not the result of reducing to the standard
form. However, replacing $g_{i-b}$ by the uniquely determined
linear combination
$$
g^+_{i-b}=g_{i-b}-\sum^{i-b-1}_{j=1}\lambda_jg_j,
$$
$dg^+_{i-b}(o)=0$, and taking the closure, we get the set
$$
\bar{C}\subset {\cal P}^{i-b-1}_{\leq 2,M-1}\times {\cal
P}^{b}_{2,M-1},
$$
which already is the result of reducing to the standard form of a
certain closed subset $C\subset {\cal P}^{i-1}_{\leq 2,M-1}$.
Taking into account the $GL_i({\mathbb C})$-invariance of the
original subvariety $B$, we conclude that all values of the
coefficients $\lambda_j$ in the formula for $g^+_{i-b}$ are
realized, so that
$$
\mathop{\rm codim} \bar{C}\leq \mathop{\rm
codim}\bar{B}_L+(i-b-1)-(M-1)\leq
$$
$$
\leq a-(M+b-i)(b+1)
$$
and for that reason
$$
\mathop{\rm codim} C\leq a-(M+b-i)(b+1)+(M-1)b-(i-b-1)b=
$$
$$
= a- (M+b-i),
$$
whereas $\varepsilon(C)=b$. Since in the case under consideration
$\alpha=0$, we get that the multiplicity of the intersection
$$
\mathop{\rm mult}\nolimits_{o}\{f_1=\dots=f_{i-1}=L=0\},
$$
as in the case 1), can be estimated from above by the number
$$
\mu_{i-1,M-1}(a-(M+b-i)-\alpha(b-1),e-\alpha),
$$
which completes the proof of Theorem 1. Q.E.D.\vspace{0.1cm}

{\bf Remark 1.4.} For $a\leq M$ the claim of Theorem 1 and its
proof remain valid for spaces of polynomials of arbitrary degree
$d\geq 2$. In the beginning of the proof of Theorem 1 (Sec. 1.4)
the polynomial $f_i$ should be taken in the form $gL_1$, where $g$
is a generic polynomial of degree $(d-1)$ (it is sufficient to
require that the differential $dg(o)$ is a linear form of general
position with respect to a generic tuple $(f_1,\dots,f_{i-1})$),
and $L_1\in{\cal P}_{1,M}$ is a linear form. The proof given above
works without any modifications.\vspace{1cm}

\newpage


{\bf \S2. Asymptotic estimates}\vspace{0.3cm}

In this section, using the inductive inequality of Theorem 1, we
obtain upper bounds for the numbers $\mu_{i,M}(a,b)$ and
$\mu_i(a)$ and consider their asymptotics for sufficiently high
values of $M$.\vspace{0.3cm}

{\bf 2.1. Estimates for the small values of $\varepsilon=b$.} As
we mentioned above, for the trivial reasons
$\mu_{i,M}(a,0)=1$.\vspace{0.1cm}

{\bf Example 2.1.} Let us obtain an upper bound for the numbers
$\mu_{i,M}(a,1)$. We get
$$
\mu_{i,M}(a,1)\leq 1+\mu_{i-1,M-1}(a-(M+1-i),1-\alpha_1).
$$
If $\alpha_1=1$, then $\mu_{i.M}(a,1)\leq 2$. If $\alpha_1=0$,
then Theorem 1 can be applied once again. Assume that the value of
the parameter $\alpha$ is 0 at the first $k$ steps:
$$
\alpha_1=\dots=\alpha_k=0.
$$
Applying Theorem 1 $k$ times, we get:
$$
\begin{array}{rccl}
\mu_{i,M}(a,1) & \leq & 1+ \mu_{i-1,M-1}(a-(M+1-i),1) & \leq \\
               & \leq & 2+ \mu_{i-2,M-2}(a-2(M+1-i),1) & \leq \\
               &      &  \dots                         &      \\
               & \leq & k+ \mu_{i-k,M-k}(a-k(M+1-i),1). &
\end{array}
$$
This is possible if the inequality
$$
a\geq (k+1)(M+1-i)
$$
holds. Therefore, the maximal possible number $k$ of steps, at
which the parameter $\alpha$ keeps the value 0, is equal to
$$
\left[\frac{a}{M+1-i}\right]-1.
$$
As a result, we obtain the estimate
$$
\mu_{i,M}(a,1)\leq \left[\frac{a}{M+1-i}\right]+1,
$$
in particular, $\mu_{M,M}(a,1)\leq a+1$. Note that the last
estimate is precise: the equality $\varepsilon=1$ means that the
complete intersection
$$
\{f_1=\dots=f_{M-1}=0\}
$$
is a smooth curve at the point $o$. The condition of tangency of
order $a\leq M$ imposes on the polynomial $f_M$ at most $a$
independent conditions. As a result we obtain the equality
$$
\mu_{M,M}(a,1)=a+1.
$$
\vspace{0.1cm}

{\bf Example 2.2.} Let us obtain an upper bound for the numbers
$\mu_{i,M}(a,2)$. Again let us assume that at the first $k$ steps
the value of the parameter $\alpha$ is equal to 0. This is
possible, if the inequality $a\geq(k+1)(M+2-i)$ holds. After $k$
applications of Theorem 1 we obtain the inequality
$$
\mu_{i,M}(a,2)\leq\sum^k_{j=1}\mu_{i-j,M-j+1}(a-j(M+2-i),1)+
$$
$$
+\mu_{i-k,M-k}(a-k(M+2-i)-1,1).
$$
Taking the maximal possible value of $k$ and using the estimate of
the previous example, we get
$$
\mu_{i,M}(a,2)\leq \frac12
\left[\frac{a}{M+2-i}\right]\left(\left[\frac{a}{M+2-i}\right]
+1\right)+2.
$$
For $i=M$ this estimate can be made slightly more precise:
$$
\mu_{M,M}(a,2)\leq \frac12\left[\frac{a}{2}\right]\left(
\left[\frac{a}{2}\right]+1\right)+\delta,
$$
where $\delta=1$, if $a$ is even, and $\delta=2$, if $a$ is odd.

In a similar way one can obtain an upper estimate for
$\mu_{i,M}(a,b)$ for $b=3,4,\dots$: applying several times Theorem
1, we can ensure that in the right hand side of the inequality the
value of the parameter $\varepsilon$ were equal to $b-1$ in all
summands, after which we can apply the inequality for
$\mu_{i,M}(a,b-1)$, obtained at the previous step.\vspace{0.3cm}

{\bf 2.2. The general method.} Applying Theorem 1 $k$ times in the
same way as we did in Examples 2.1 and 2.2, under the assumption
that the value of the parameter $\alpha$ is equal to 0, we obtain
the inequality
$$
\mu_{i,M}(a,b)\leq\sum^k_{j=1}\mu_{i-j,M-j+1}(a-j(M+b-i),b-1)+
$$
$$
+\mu_{i-k,M-k}(a-k(M+b-i)-(b-1),b-1).
$$
Note that the inequality $a\geq(k+1)(M+b-i)$ holds. However, it is
difficult to obtain in this way a general estimate for
$\mu_{i,M}(a,b)$, reducing it to the estimate for the numbers with
$\varepsilon=b-1$, because of the difficult formulas, which are
hard to follow. However, we may conclude that a multiple
application of Theorem 1 yields the estimate
\begin{equation}\label{june2011.3}
\mu_{i,M}(a,b)\leq\sum_{j,N,a',b'}\mu_{j,N}(a',b')
\end{equation}
for a certain set of tuples $(j,N,a',b')$ (possibly, with
repetitions of the same tuple), and in the end, the estimate
\begin{equation}\label{june2011.4}
\mu_{i,M}(a,b)\leq\sum_{j,N,a'}\mu_{j,N}(a',0),
\end{equation}
where in the right hand side all components are equal to 1, so
that it is sufficient to estimate from above the number of
components, which is equal to the number of inductive steps ---
applications of Theorem 1. For this purpose, with each term in the
right hand side of the inequality (\ref{june2011.3}) we associate
a word
$$
\omega=\tau_1\tau_2\dots\tau_K
$$
in the alphabet $\{A,B_0,B_1\}\ni \tau_i$, describing the
``origin'' of that term. With the term $\mu_{i,M}(a,b)$ itself in
the tautological estimate
$$
\mu_{i,M}(a,b)\leq\mu_{i,M}(a,b)
$$
we associate the empty word. Let
\begin{equation}\label{june2011.5}
\mu_{i,M}(a,b)\leq\sum_{w\in W'}\mu[w]
\end{equation}
be the new writing of the inequality (\ref{june2011.3}), where
each term $\mu_{j,N}(a',b')$ in the right hand side corresponds to
a word $w\in W'$ and is written as $\mu[w]$. Let us choose and fix
such a term with $b'\geq 1$. According to the proof of Theorem 1,
this term gives an upper estimate for the number $\mu(B')$, where
$B'\subset {\cal P}^{j}_{\leq 2,N}$ a certain $GL_j({\mathbb
C})$-invariant irreducible subvariety of codimension $a'$ with
$\varepsilon(B')=b'$. Now, applying Theorem 1, we replace (keeping
the inequality) the term $\mu_{j,N}(a',b')$ by the sum of two new
numbers $\mu[w_1]+\mu[w_2]$, where $\mu[w_1]$ and $\mu[w_2]$
correspond to the first and second terms in the right hand side of
the inequality (\ref{june2011.1}), respectively. Here $w_1=wA$ and
$w_2=wB_{\alpha}$, where $\alpha=\alpha(B)\in\{0,1\}$. This
determines the procedure of constructing the words $w$ in a unique
way. It is clear that with each word at most one term in
(\ref{june2011.3}) is associated. Thus to obtain an upper estimate
for $\mu_{i,M}(a,b)$, we need to estimate the number of words, to
which terms in the inequality (\ref{june2011.4}) correspond.

For instance, in Example 2.1 the set of words is
$$
A,\,\,B_0A,\,\,\dots,\,\,  \underbrace{B_0\dots B_0}_k A,\,\,
\underbrace{B_0\dots B_0}_k B_1.
$$
\vspace{0.1cm}

{\bf Remark 2.1.} Let $\nu\colon\{A,B_0,B_1\}\to\{A,B\}$ be the
map of the three-letter alphabet into the two-letter one, given by
$\nu(A)=A$, $\nu(B_{\alpha})=B$,
$$
\nu\colon w=\tau_1\dots\tau_K\mapsto {\bar
w}=\nu(\tau_1)\dots\nu(\tau_K)
$$
the corresponding map of the set of words. Then for any inequality
(\ref{june2011.5}), obtained by an application of Theorem 1, the
restriction $\nu|_{W'}$ is injective. Indeed, each application of
Theorem 1 replaces some word $w$ by the pair of words $wA$ and
$wB_{\alpha}$, where the value of the parameter $\alpha$ is
uniquely determined.\vspace{0.1cm}

Now with each summand $\mu_{j,N}(a',b')$ (or with the word $w$,
corresponding to that summand) we associate the triple of
integer-valued parameters $(a',b',\Delta')$, where
$\Delta'=N+b'-j$. By Theorem 1,

\begin{itemize}

\item for the word $wA$ the associated triple is $(a'-\Delta',b'-1,\Delta')$,

\item for the word $wB_0$ it is the triple $(a'-\Delta',b',\Delta')$,

\item for the word $wB_1$ it is the triple $(a'-\Delta'-(b'-1),b'-1,\Delta'-1)$.

\end{itemize}

Recall now that the term $\mu_{j,N}(a',b')$ is well defined only
if the inequality $a'\geq b'\Delta'$ holds.

Let $W$ be the set of words, corresponding to the summands of the
right hand side of the inequality (\ref{june2011.4}). Let
$W_l\subset W$ be the subset, consisting of the words, in which
precisely $l$ letters are $B_1$. Obviously,
$$
W=\coprod^b_{l=0}W_l
$$
(the union is disjoint), so that
$$
\sharp W=\sum^b_{l=0}\sharp W_l.
$$
It remains to estimate from above the number of elements in each
of the sets $W_l$.\vspace{0.1cm}

{\bf Lemma 2.1.} {\it The inequality
$$
\sharp W_l\leq\left(\begin{array}{c}
A_l\\
b-l \end{array}\right)
$$
holds, where} $\displaystyle
A_l=\left[\frac{a-lb}{\Delta-l}\right]$.\vspace{0.1cm}

{\bf Proof.} Consider first the case $l=0$. In the word $w\in W_l$
there are no letters $B_1$, whereas the letter $A$ occurs
precisely $b$ times, since to the word $w$ corresponds the triple
$(a',0,\Delta')$, and the letter $B_0$ does not change the value
of the parameter $\varepsilon=b'$. On the other hand, since the
letter $B_1$ does not occur, we get $\Delta'=\Delta=M+b-i$, and
the inequality $a'\geq 0$ implies that the length of the word $w$
does not exceed $A_0=[a/\Delta]$. Thus $\sharp W_0$ does not
exceed the number of ways of putting $b$ letters $A$ on at most
$A_0$ positions. However, the last letter in the word $w\in W_0$
can be only the letter $A$, by the same reason that $A$ decreases
the value of $\varepsilon=b'$ by 1, and $B_0$ does not change it.
Therefore, $\sharp W_0$ does not exceed the number of ways of
putting $b$ letters $A$ on $A_0$ positions, which is what we
need.\vspace{0.1cm}

Now let us consider the case of an arbitrary $l\leq
b$.\vspace{0.1cm}

{\bf Lemma 2.2.} {\it The length of a word $w\in W_l$ does not
exceed $A_l$}.\vspace{0.1cm}

Accepting the claim of Lemma 2.2, let us complete the proof of
Lemma 2.1. Obviously, the letter $A$ occurs in a word $w\in W_l$
precisely $(b-l)$ times. We associate with the word $w$ the
corresponding way of putting $(b-l)$ letters $A$ on $A_l$
positions.

We claim that this map is injective. (This immediately implies
Lemma 2.1.) Indeed, assume that this is not true: there are two
distinct words $w_1\neq w_2$ in $W_l$ with the same distribution
of the letter $A$. Assume that the length $|w_1|$ of the word
$w_1$ does not exceed the length $|w_2|$. Changing to the
two-letter alphabet $\{A,B\}$, we conclude that the letter $w_1$
is a left segment of the word $w_2$ and
$$
w_2=w_1B_{\alpha}\dots B_{\alpha_{k}}
$$
for some $\alpha_1,\dots,\alpha_k\in\{0,1\}$. However, the
parameter $\varepsilon=b'$ of the word $w_1$ is already equal to
0, which implies that $w_1=w_2$. Q.E.D. for Lemma
2.1.\vspace{0.1cm}

{\bf Proof of Lemma 2.2.} Let us control the length $|w|$ of the
word $w\in W_l$ by the decreasing of the parameter $a'\geq 0$. The
slower it decreases, the longer can be the word. Assume that the
letter $B_1$ occupies the positions
$$
k_1,k_1+k_2,\dots,k_1+k_2+\dots+k_l,
$$
where $k_i\geq 1$. On each segment
$$
[k_1+\dots +k_j+ 1,k_1+\dots+k_{j+1}-1]
$$
of the word $w$ (provided it is non-empty) the value of the
parameter $\varepsilon=b'$ can get smaller by, at most,
$k_{j+1}-1$, whereas the value of the parameter $\Delta'$ remains
the same. Therefore, to the left segment of the word $w$ of length
$k_1+\dots+k_l$ corresponds the value
$$
\begin{array}{cll}
a'\geq a & -k_1\Delta-(b-k_1)-\\
& -k_2(\Delta-1)-(b-k_1-k_2)-\\
& \dots\\
& -k_l(\Delta-(l-1))-(b-k_1-\dots-k_l)=\\
& =(a-lb)-(\Delta-l)(k_1+\dots+k_l).
\end{array}
$$
After the position $(k_1+\dots+k_l)$ the value of the parameter
$\Delta'$ remains the same and is equal to $(\Delta-l)$.
Therefore,
$$
|w|\leq k_1+\dots+k_l+
\left[\frac{(a-lb)-(\Delta-l)(k_1+\dots+k_l)}{\Delta-l}\right]=A_l.
$$
Q.E.D. for Lemma 2.2.\vspace{0.1cm}

{\bf Corollary 2.1.} {\it The inequality
$$
\mu_{i,M}(a,b)\leq\sum^b_{l=0}\left(\begin{array}{c}A_l\\
b-l\end{array}\right)
$$
holds, where $\displaystyle
A_l=\left[\frac{a-lb}{M+b-i-l}\right]$}.\vspace{0.3cm}

{\bf  2.3. An asymptotic estimate for a high dimension.} Obtaining
compact upper estimates for the numbers $\mu_i(a)$, which could be
used for particular computations, presents a non-trivial problem.
The inequality of Corollary 2.1 is too complicated and not very
visual. However, in one case it is easy to derive from it a simple
and precise estimate.\vspace{0.1cm}

{\bf Example 2.3.} Assume that $M=m^2$ is a full square. Then the
following equality holds:
$$
\mu_{M,M}(M,m)=2^m.
$$
Indeed, all numbers $A_l=m$ are the same, so that we get
$$
\mu_{M,M}(M,m)\leq\sum^m_{l=0}\left(\begin{array} {c}m\\
l\end{array}\right)=2^m.
$$
On the other hand, obviously $\mu_{M,M}(M,m)\geq 2^m$.
Q.E.D.\vspace{0.1cm}

Now let us consider the general case for $i=M$ and the maximal
possible codimension $a=M$. Set $\xi(M)=\mu_M(M)$. Since
$$
\xi(M)=\mathop{\rm max}\limits_{1\leq b\leq
[\sqrt{M}]}\mu_{M,M}(M,b),
$$
by Corollary 2.1 we get
$$
\xi(M)\leq\sqrt{M}\max\left(\begin{array}{c}\displaystyle
\left[\frac{M-lb}{b-l}\right]
\\ \\b-l
\end{array}\right),
$$
where the maximum is taken over $b\in\{1,\dots,[\sqrt{M}]\}$ and
$l\in\{1,\dots,b\}$. Now elementary computations with binomial
coefficients and an application of the Stirling formula give the
following result. Set
$$
\omega=\mathop{\rm max}\limits_{s\in[1,\infty)}[2s\mathop{\rm
ln}s-(s-\frac{1}{s})\mathop{\rm ln}(s^2-1)].
$$
\vspace{0.1cm}

{\bf Proposition 2.1.} {\it For sufficiently high $M$ the
inequality
$$
\xi(M)\leq\sqrt{M}e^{\omega\sqrt{M}}
$$
holds, where $e$ is the base of the natural
logarithm.}\vspace{1cm}

\newpage


\begin{center}{\bf \S3. Systems of equations with the set of solutions
\\ of ``incorrect'' dimension}
\end{center}\vspace{0.3cm}

In this section, we prove Proposition 1.1.\vspace{0.3cm}

{\bf 3.1. Systems of homogeneous equations.} In the space ${\cal
P}^{i}_{d,M+1}$ of systems of homogeneous polynomials
$(p_1,\dots,p_i)$ of degree $d\geq 2$ in the variables
$z_0,\dots,z_M$ consider the closed subset $Y$, consisting of such
tuples $(p_1,\dots,p_i)$, that the set
$$
\{p_1=\dots=p_i=0\}\subset {\mathbb P}^M
$$
is of ``incorrect'' codimension $\leq i-1$.\vspace{0.1cm}

{\bf Proposition 3.1.} {\it The codimension of the subset $Y$ in
the space ${\cal P}^{i}_{d,M+1}$ is not less than}
$$
\mathop{\rm
min}\limits_{b\in\{0,\dots,i-1\}}\{((b+1)d-b)(M-b)+1\}.
$$
\vspace{0.1cm}

{\bf Proof.} It follows directly from [3, Proposition 4], taking
into account that the degrees of the polynomials $(p_1,\dots,p_i)$
are equal. Q.E.D.\vspace{0.1cm}

{\bf Corollary 3.1.} {\it For $i\leq M-1$ the codimension of the
subset $Y$ in the space ${\cal P}^{i}_{d,M+1}$ is not less than
$dM+1$, and for $i=M$ it is not less than
$(d-1)M+2$}.\vspace{0.1cm}

{\bf Proof.} Since $d\geq 2$, the quadratic function
$$
\gamma(b)=b^2(1-d)+b(dM-M-d)+dM+1
$$
of the variable $b$ is negative definite and attains its maximum
at
$$
b_*=\frac{dM-M-d}{2(d-1)}>0.
$$
Therefore, the minimum of this function on the set
$\{0,\dots,i-1\}$ is attained either for $b=0$ (and equal to
$dM+1$), or for $b=i-1$. It is easy to check that
$\gamma(M-1)=2(dM-d-M)+5\geq dM+1$, which proves the first claim
of the corollary. Furthermore, $\gamma(M)=(d-1)M+2\leq dM+1$,
which proves the second claim. Q.E.D.\vspace{0.3cm}

{\bf 3.2. Systems of non-homogeneous equations.} Let us prove
Proposition 1.1. In the space ${\cal P}^{i}_{d,M+1}\times {\mathbb
P}^M$ consider the closed algebraic set ${\cal Y}$, consisting of
such pairs $((p_1,\dots,p_i), x\in{\mathbb P}^M)$, that the
corresponding set of zeros $\{p_1=\dots=p_i=0\}$ has an
irreducible component of ``incorrect'' codimension $\leq i-1$,
passing through the point $x$. Furthermore, denote by the symbol
$Y_b$, $b=0,\dots,i-1$, the closed subset in $Y$, consisting of
such tuples $(p_1,\dots,p_i)$, that the codimension of the set of
zeros $\{p_1=\dots=p_i=0\}$ does not exceed $i-1-b$; in
particular, $Y_0=Y$. By the methods of [3, Sec. 3] it is easy to
check that $\mathop{\rm codim}_YY_b\geq 2b$ (in fact, the estimate
is much stronger). This implies that
$$
\mathop{\rm dim{\cal Y}}=\mathop{\rm dim}Y+M-i+1.
$$
By the symbols $\pi_1$ and $\pi_2$ denote the projections of the
direct product ${\cal P}^{i}_{d,M+1}\times {\mathbb P}^M$ onto the
first and second factors, respectively. Obviously, $\pi_1({\cal
Y})=Y$. Furthermore, $\pi_2({\cal Y})={\mathbb P}^M$, and all the
fibres $\pi^{-1}_2(x)\cap{\cal Y}={\cal Y}_x$ are of the same
dimension
$$
\mathop{\rm dim}{\cal Y}-M=\mathop{\rm dim} Y-i+1.
$$
On the other hand, the space ${\cal P}^{i}_{\leq d,M}$ can be
naturally identified with the closed subset of codimension $i$ in
$\pi^{-1}_2(x)\cong{\cal P}^i_{d,M+1}$, consisting of such tuples
$(p_1,\dots,p_i)$, that
$$
p_1(x)=\dots=p_i(x)=0.
$$
It is clear that ${\cal Y}_x$ is contained in that subset, so that
the codimension of ${\cal Y}_x$ with respect to ${\cal P}^i_{\leq
d,M}$ is equal to
$$
\mathop{\rm codim}(Y\subset {\cal P}^{i}_{d,M+1})-1.
$$
Applying Corollary 3.1, we complete the proof of Proposition
1.1.\vspace{0.3cm}

{\bf 3.3. Precision of the estimates.} How precise are the
estimates of Proposition 1.1? The following example shows that for
$i=M$ the estimate is sharp. Let $L\ni o$ be an arbitrary line
passing through the origin. The condition that
$$
p(z_1,\dots,z_M)|_L\equiv 0
$$
imposes on a polynomial of degree $d$ precisely $d$ independent
conditions (recall that $p(0,\dots,0)=0)$. Therefore, requiring
that
$$
L\subset\{p_1=\dots=p_M=0\}
$$
we impose on the tuple of polynomials $(p_1,\dots,p_M)\in {\cal
P}^{M}_{\leq d,M}$ precisely $dM$ independent conditions. Since
there is a $(M-1)$-dimensional family of lines, passing through
the point $o$, the set of tuples $(p_1,\dots,p_M)$ such that the
closed set $\{(p_1=\dots=p_M=0)\}$ contains a line passing through
the point $o$, is of codimension $(d-1)M+1$ in the space ${\cal
P}^{M}_{\leq d,M}$. Therefore, the estimate of Proposition 1.1 is
sharp. In particular, the set of tuples $(p_1,\dots,p_M)$,
vanishing on a line, forms an irreducible component of the set
$Y$. The question, what is the codimension of other components of
this set, remains an open problem.\vspace{1cm}


{\bf References}\vspace{0.3cm}

{\small

\noindent 1. Pukhlikov A.V., Birationally rigid varieties with a
pencil of Fano double covers. II,  Sbornik: Mathematics. V. 195
(2004), no. 11, 1665-1702. \vspace{0.1cm}

\noindent 2. Fulton W., Intersection Theory, Springer-Verlag,
1984.\vspace{0.1cm}

\noindent 3. Pukhlikov A.V., Birationally rigid Fano complete
intersections, Crelle J. f\" ur die reine und angew. Math. {\bf
541} (2001), 55-79. \vspace{0.1cm}

}

\begin{flushleft}
\it pukh@liv.ac.uk \\
\it pukh@mi.ras.ru
\end{flushleft}

\end{document}